\documentclass[11 pt,a4paper,leqno]{article}
\usepackage{latexsym,amsmath,amscd,amsfonts}

\RequirePackage{amsmath}
\RequirePackage{amssymb}
\RequirePackage{ifthen}
\RequirePackage{theorem}
\RequirePackage{pslatex}

\usepackage[all]{xy}

\title{\sc Classification de modules aux diff\'erences filtr\'es isogradu\'es}

\author{J. Sauloy
\footnote{Laboratoire \'Emile Picard, 
Institut de Math\'ematiques, CNRS UMR 5219, 
U.F.R. M.I.G., Universit\'e Paul Sabatier (Toulouse 3),
31062 Toulouse CEDEX 9}}

\date{}

\def\C{{\mathbf C}}

\def\Z{{\mathbf Z}}

\def\N{{\mathbf N}}

\def\F{{\mathcal{F}}}

\def\B{{\mathcal{B}}}

\def\D{{\mathcal{D}}}

\def\gr{{\text{gr}}}

\def\G{{\mathfrak{G}}}

\def\Raq(d){{\C\{\xi\}_{q,(\delta)}}}
\def\Kaq(d){{\C(\{\xi\})_{q,(\delta)}}}

\def\DM{{DiffMod(K,\sigma)}}

\def\Id{\text{Id}}

\def\Aut{\text{Aut~}}

\newtheorem{thm}{Theor\`eme}[section]
\newtheorem{lem}[thm]{Lemme}
\newtheorem{prop}[thm]{Proposition}
\newtheorem{cor}[thm]{Corollaire}
\newtheorem{defn}[thm]{D\'efinition}

\def\Pr{\textsl{Preuve. - }}
\def\Endpr{$\Box$ \hfill \\}
\def\Ex{\noindent\textbf{Exemple.~}}

\def\Rem{\noindent\textbf{Remarque.~}}

\def\tq{\, | \,}
\def\Lin{\mathcal{L}}
\def\Hom{\text{Hom}}
\def\Ext{\text{Ext}}
\def\1{\underline{1}}
\def\ie{\emph{i.e.}}

\begin{document}

\maketitle

\bigskip \hrule \bigskip

\begin{abstract}
{\small 
La classification analytique locale des \'equations aux $q$-diff\'erences
irr\'eguli\`eres (\cite{RSZ}) se ram\`ene \`a la classification de modules
aux $q$-diff\'erences filtr\'es \`a gradu\'e fix\'e. Nous d\'egageons ici
des hypoth\`eses g\'en\'erales qui assurent l'existence d'un sch\'ema de
modules pour ce probl\`eme, qui soit de plus un espace affine.
}
\end{abstract}

\bigskip \hrule \bigskip

%%%%%%%%%%%%%%%%%%%%%%%%%%%%%%%%%%%%%%%%%%%%%%%%%%%%%%%%%%%%%%%%%%%%%%%%%%%%%

\section*{Formulation g\'en\'erale du probl\`eme}

Les \emph{modules aux $q$-diff\'erences} associ\'es aux \'equations aux 
$q$-diff\'erences lin\'eaires analytiques en $0$ sont naturellement munis 
d'une filtration par les pentes respect\'ee par les morphismes analytiques
et telle que deux modules sont formellement isomorphes si, et seulement
s'ils ont des gradu\'es isomorphes. La \emph{classification analytique
isoformelle} (\cite{RSZ}) est ainsi ramen\'ee \`a la \emph{classification 
isogradu\'ee de modules filtr\'es}. Cette derni\`ere admet une formulation
purement alg\'ebrique et pose un probl\`eme assez naturel, que nous allons
d\'ecrire sous une forme un peu plus g\'en\'erale. \\

Soient $C$ un anneau commutatif et $\mathcal{C}$ une cat\'egorie
ab\'elienne $C$-lin\'eaire. On fixe un objet finiment gradu\'e:
$$
P = P_{1} \oplus \cdots \oplus P_{k},
$$
et l'on se propose de classifier les couples 
$(\underline{M},\underline{u})$ form\'es d'un objet finiment filtr\'e: 
$$
\underline{M} = 
(0 = M_{0} \subset M_{1} \subset \cdots \subset M_{k} = M),
$$
et d'un isomorphisme de $\gr M$ sur $P$:
$$
\underline{u} = (u_{i}: M_{i}/M_{i-1} \simeq P_{i})_{1 \leq i \leq k}.
$$
Il revient d'ailleurs au m\^eme de se donner $k$ suites exactes:
$$
0 \rightarrow M_{i-1} \rightarrow M_{i} 
\overset{v_{i}}{\rightarrow} P_{i} \rightarrow 0.
$$
(On laisse au lecteur le soin de se convaincre que cela revient en effet
au m\^eme.) Les couples $(\underline{M},\underline{u})$ et
$(\underline{M'},\underline{u'})$ sont d\'eclar\'es \'equivalents
si (avec les notations \'evidentes pour $\underline{M'}$) il existe 
un morphisme de $M$ dans $M'$ qui respecte les filtrations et 
compatible avec les isomorphismes des gradu\'es, \ie\ tel que 
$u = u' \circ \gr f$, \ie\ tel que le diagramme ci-dessous soit
commutatif:
$$
\xymatrix{
\gr M \ar@<0ex>[rd]^{u} \ar@<0ex>[rr]^{\gr f}  &  
& \gr M'  \ar@<0ex>[ld]_{u'} \\
& P &  \\
}
$$
Dans la description par suites exactes: s'il 
existe des morphismes $f_{i}: M_{i} \rightarrow M'_{i}$ rendant
commutatifs les diagrammes:
$$
\begin{CD}
0 @>>> M_{i-1}   @>>> M_{i}     @>{v_{i}}>>      P_{i}  @>>> 0     \\
&  &   @VV{f_{i-1}}V      @VV{f_{i}}V      @VV{\Id_{P_{i}}}V       \\
0 @>>> M'_{i-1}   @>>> M'_{i}     @>{v'_{i}}>>    P_{i} @>>> 0
\end{CD}
$$

\Rem
S'il existe un tel morphisme, il est automatiquement strict et
un isomorphisme. \\

On notera $\F(P_{1},\ldots,P_{k})$ l'ensemble des classes de couples
$(\underline{M},\underline{u})$ pour la relation d'\'equivalence
ci-dessus. Nous admettons que c'est bien un ensemble: on se restreint 
si l'on y tient \`a de petites cat\'egories, ou bien on suppose que les
$\Ext$ en sont, ce qui (d'apr\`es les arguments de d\'evissage qui vont 
suivre page \pageref{methode}) suffit. Dans le cas que nous traiterons 
(sous-cat\'egories de la cat\'egorie des modules \`a gauche sur un anneau), 
c'est facile \`a v\'erifier. \\

Pour $k = 1$, l'ensemble $\F(P_{1})$ est un singleton. Pour $k = 2$, 
l'ensemble $\F(P_{1},P_{2})$ s'identifie naturellement \`a l'ensemble 
$\Ext(P_{2},P_{1})$ des classes d'extensions de $P_{2}$ par $P_{1}$,
lequel porte une structure de $C$-module. L'identification s'obtient
comme suit. Se donner un module filtr\'e  
$\underline{M} = (0 = M_{0} \subset M_{1} \subset M_{2} = M)$ muni
d'un isomorphisme de $\gr M$ sur $P_{1} \oplus P_{2}$ revient \`a
se donner un isomorphisme de $M_{1}$ sur $P_{1}$ et un isomorphisme de 
$M/M_{1}$ sur $P_{2}$, soit encore un monomorphisme $i$ de $P_{1}$
dans $M$ et un \'epimorphisme $p$ de $M$ sur $P_{2}$ de noyau
$i(P_{1})$, soit encore une suite exacte
$$
0 \rightarrow P_{1} \overset{i}{\rightarrow} M 
\overset{p}{\rightarrow} P_{2} \rightarrow 0,
$$
c'est-\`a-dire une extension de $P_{2}$ par $P_{1}$. Et l'on v\'erifie
sans peine que notre relation d'\'equivalence correspond ainsi \`a
l'isomorphisme usuel d'extensions. \\

\label{methode}
Notre but est de trouver des conditions qui garantissent que 
$\F(P_{1},\ldots,P_{k})$ porte une structure d'espace affine sur $C$,
et d'en calculer la dimension. Le cas $k = 2$ sugg\`ere de supposer
que les $C$-modules $\Ext(P_{j},P_{i})$ sont libres de rang fini
(on verra qu'en fait seuls comptent les couples tels que $i < j$).
On proc\`ede ensuite \`a un d\'evissage. En vue d'une r\'ecurrence 
sur $k$, on invoque une surjection naturelle:
$$
\F(P_{1},\ldots,P_{k}) \longrightarrow \F(P_{1},\ldots,P_{k-1}),
$$
qui, \`a la classe de $(\underline{M},\underline{u})$ d\'efini
comme plus haut, associe la classe de $(\underline{M'},\underline{u'})$
d\'efini par 
$$
\underline{M'} = 
(0 = M_{0} \subset M_{1} \subset \cdots \subset M_{k-1} = M')
\text{~et~}
\underline{u'} = (u_{i})_{1 \leq i \leq k-1}.
$$
L'image r\'eciproque de la classe de $(\underline{M'},\underline{u'})$
d\'ecrit comme ci-dessus s'identifie \`a $\Ext(P_{k},M')$. Noter que 
$\Ext(P_{k},M')$ ne d\'epend bien en effet (\`a isomorphisme canonique 
pr\`es) que de la classe de $(\underline{M'},\underline{u'})$. Sous
les conditions o\`u nous nous placerons, on verra que $\Ext(P_{k},M')$
se d\'evisse \`a son tour en les $\Ext(P_{k},P_{i})$ pour $i < k$, et
l'on s'attend \`a obtenir un espace de dimension
$\sum\limits_{1 \leq i < j \leq k} \dim \Ext(P_{j},P_{i})$. \\

\Rem
Une fois d\'ecrit l'espace $\F(P_{1},\ldots,P_{k})$, on peut poser le 
probl\`eme (apparemment plus naturel) de la classification des $M$ tels 
que $\gr M \simeq P$ (o\`u l'on ne prescrit pas la ``polarisation'' $u$). 
On v\'erifie que $\prod \Aut(P_{i})$ agit sur $\F(P_{1},\ldots,P_{k})$
(c'est l'action sur ``l'ensemble'' des $(\underline{M},\underline{u})$
par compositions \`a gauche $\phi_{i} \circ u_{i}$) et qu'il s'agit de 
quotienter par cette action. On ne s'occupera pas ici de ce probl\`eme. \\

L'usage de l'alg\`ebre homologique dans les probl\`emes de classification 
des \'equations diff\'erentielles n'est pas nouveau, mais l'auteur de ce
texte croit avoir remarqu\'e que l'\'etape de mod\'elisation alg\'ebrique
pr\'eliminaire est souvent trait\'ee de mani\`ere d\'esinvolte. Ainsi,
l'identification des modules d'extensions \`a des conoyaux comporte
presque toujours la description explicite d'une application, parfois la 
preuve de sa bijectivit\'e, rarement la preuve de son additivit\'e et
(semble-t-il) jamais la preuve de sa lin\'earit\'e. Pour cette raison,
un tr\`es grand soin a \'et\'e ici apport\'e au d\'etail des constructions
alg\'ebriques et aux d\'emonstrations d'isomorphies "\'evidentes".

\paragraph{Remerciements.}
Ce texte a \'et\'e achev\'e lors d'un s\'ejour en tant que professeur
invit\'e au laboratoire "Algebra, Geometria y Topologia" de l'Universit\'e
de Valladolid, avec le soutien de l'Instituto de Estudios Iberoamericanos.
Je remercie chaleureusement Jose-Manuel Aroca, Felipe Cano y los otros
amigos de ahi pour leur accueil amical.

% 1 
 
\section{Modules aux diff\'erences}
\label{section:moddiff}

Soient $C$ un anneau commutatif
\footnote{Le choix de la lettre $C$ peut \'evoquer le corps $\C$
des complexes ou bien un corps de ``constantes''.}, 
$K$ une $C$-alg\`ebre commutative et $\sigma$ un automorphisme de 
la $C$-alg\`ebre $K$. On note:
$$
K^{\sigma} := \{x \in K \tq \sigma x = x\}
$$
la sous-$C$-alg\`ebre des \emph{constantes}, et l'on introduit
l'anneau de \"Ore des \emph{op\'erateurs aux diff\'erences}:
$$
\D_{K,\sigma} := K<\sigma,\sigma^{-1}> = 
\Bigl\{\sum_{i \in \Z} a_{i} \sigma^{i} \tq (a_{i}) \in K ^{(\Z)}\Bigr\}.
$$
Ses \'el\'ements sont des polyn\^omes de Laurent non commutatifs,
et la multiplication est caract\'eris\'ee par la relation de
commutation tordue:
$$
\forall k \in \Z \;,\; \lambda \in K \;,\; 
\sigma^{k}.\lambda = \sigma^{k}(\lambda) \sigma^{k}.
$$
Le centre de la $C$-alg\`ebre $\D_{K,\sigma}$ est $K^{\sigma}$. \\

\Rem
En principe, la th\'eorie des modules aux diff\'erences repose 
sur celle des corps diff\'erentiels: $K$ serait un corps et 
$C := K^{\sigma}$ son corps des constantes. Nos hypoth\`eses
plus faibles sont choisies en vue de la situation "relative"
(extension des scalaires).

\begin{defn}
Un \emph{module aux diff\'erences} sur la \emph{$C$-alg\`ebre 
aux diff\'erences} $(K,\sigma)$ (ou, pour faire court, sur $K$)
est un $\D_{K,\sigma}$-module \`a gauche qui induit, par 
restriction des scalaires, un $K$-module projectif de rang fini.
La sous-cat\'egorie pleine de $\D_{K,\sigma}-Mod$ (cat\'egorie
des $\D_{K,\sigma}$-modules \`a gauche) form\'ee des modules 
aux diff\'erences sera not\'ee $\DM$. (Les deux cat\'egories
sont ab\'eliennes et $C$-lin\'eaires.)
\end{defn}

Tout $\D_{K,\sigma}$-module \`a gauche peut \^etre r\'ealis\'e
par un couple $(E,\Phi)$, o\`u $E$ est un $K$-module (obtenu
par restriction des scalaires) et $\Phi$, qui incarne la
multiplication \`a gauche par $\sigma$, est un \emph{automorphisme
semi-lin\'eaire} ou encore \emph{$\sigma$-lin\'eaire}, \ie\ un
automorphisme du groupe $E$ tel que:
$$
\forall \lambda \in K \;,\; \forall x \in E \;,\;
\Phi(\lambda x) = \sigma(\lambda) \Phi(x).
$$
Dans cette description, un morphisme de $\D_{K,\sigma}$-modules 
\`a gauche de $(E,\Phi)$ dans $(F,\Psi)$ est une application
$u \in \Lin_{K}(E,F)$ (autrement dit, une application $K$-lin\'eaire 
$u: E \rightarrow F$) telle que $\Psi \circ u = u \circ \Phi$. 
Lorsque $E$ est projectif de rang fini, nous noterons $M := (E,\Phi)$
le module aux diff\'erences correspondant. \\

\Ex
Soient $n \in \N^{*}$ et $A \in GL_{n}(K)$. L'application
$$
X \mapsto \Phi_{A}(X) := A^{-1} \, \sigma X
$$
de $E := K^{n}$ dans lui-m\^eme est un automorphisme semi-lin\'eaire
et $M := (K^{n},\Phi_{A})$ est un module aux diff\'erences sur $K$.
Les ``vecteurs horizontaux'' de $M$ sont les $X \in K^{n}$ tels que 
$\Phi_{A}(X) = X$, autrement dit, les solutions dans $K^{n}$ de 
``l'\'equation aux $\sigma$-diff\'erences'' $\sigma X = A X$: c'est 
pour cela que nous avons utilis\'e $A^{-1}$.

% 1.1 

\subsection{Description matricielle des modules aux diff\'erences}
\label{subsection:descmat1}

On suppose ici que le $K$-module $E$ est \emph{libre} de rang fini $n$.
Cette hypoth\`ese sera reprise aux paragraphes \ref{subsection:descmat2}
et\ref{subsection:descmat3}, o\`u la description matricielle sera
poursuivie. \\

Le choix d'une base $\B$ permet d'identifier $E$ \`a $K^{n}$. Il est alors
clair que $\Phi(\B)$ est encore une base de $E$, d'o\`u l'existence de
$A \in GL_{n}(K)$ telle que $\Phi(\B) = \B A^{-1}$. Si le vecteur colonne 
des coordonn\'ees de $x \in E$ dans la base $\B$ est $X \in K^{n}$, ce
que l'on peut \'ecrire $x = \B X$, le calcul:
$$
\Phi(x) = \Phi(\B X) = \Phi(\B) \sigma(X) = \B A^{-1} \sigma(X)
$$
montre que le vecteur colonne des coordonn\'ees de $\Phi(x)$ dans $\B$ 
est $A^{-1} \sigma(X)$; autrement dit, avec les notations de l'exemple
pr\'ec\'edent, on obtient une identification:
$$
(E,\Phi) \simeq (K^{n},\Phi_{A}).
$$
Notons enfin que les morphismes de $(K^{n},\Phi_{A})$ dans $(K^{p},\Phi_{B})$ 
s'identifient aux matrices $F \in M_{p,n}(K)$ telles que:
$$
(\sigma F) A = B F,
$$ 
et la composition se ram\`ene au produit matriciel. En particulier:
$$
(K^{n},\Phi_{A}) \simeq (K^{p},\Phi_{B}) \Longleftrightarrow
n = p \text{~et~} \exists F \in GL_{n}(K) \;:\;
B = F[A] := (\sigma F) A F^{-1};
$$
(``Transformation de jauge''.)

% 1.2

\subsection{Description matricielle des modules aux diff\'erences filtr\'es}
\label{subsection:descmat2}

Soient $P_{i} = (G_{i},\Psi_{i})$ ($1 \leq i \leq k$) des module aux 
diff\'erences tels que les $K$-modules $G_{i}$ soient libres de rang
fini $r_{i}$ . On munit chaque $G_{i}$ d'une base $\mathcal{D}_{i}$, 
et l'on note $B_{i} \in GL_{r_{i}}(K)$ la matrice inversible telle que
$\Psi_{i}(\mathcal{D}_{i}) = \mathcal{D}_{i} B_{i}$. \\

Soit $\underline{M}$ un module aux diff\'erences filtr\'e de gradu\'e 
$P = P_{1} \oplus \cdots \oplus P_{k}$; plus pr\'ecis\'ement,
$\underline{M} = (0 = M_{0} \subset M_{1} \subset \cdots \subset M_{k} = M)$
est muni d'un isomorphisme 
$\underline{u} = (u_{i}: M_{i}/M_{i-1} \simeq P_{i})_{1 \leq i \leq k}$
de $\gr M$ sur $P$. Notant $M_{i} = (E_{i},\Phi_{i})$, on construit par
r\'ecurrence sur $i = 1,\ldots,k$ une base $\mathcal{B}_{i}$ de $E_{i}$
telle que $\mathcal{B}_{i-1} \subset \mathcal{B}_{i}$ et que 
$\mathcal{B}'_{i} := \mathcal{B}_{i} \setminus \mathcal{B}_{i-1}$ rel\`eve 
$\mathcal{D}_{i}$ via $u_{i}$ (on voit au passage que les $K$-modules $E_{i}$ 
soient libres de rang fini $r_{1} + \cdots + r_{i}$). \\

Notant $\Phi_{i}$ (resp. $\overline{\Phi_{i}}$) l'automorphisme 
semi-lin\'eaire induit par $\Phi$ (resp. par $\Phi_{i}$) sur $E_{i}$ 
(resp. sur $E_{i}/E_{i-1}$), et $\mathcal{C}_{i}$ la base de $E_{i}/E_{i-1}$ 
induite par $\mathcal{B}'_{i}$, on d\'eduit de l'\'egalit\'e
$u_{i} \circ \overline{\Phi_{i}} = \Psi_{i} \circ u_{i}$ (d\^ue au fait que
$u_{i}$ est un morphisme) la relation:
$\overline{\Phi_{i}}(\mathcal{C}_{i}) = \mathcal{C}_{i} B_{i}$, puis, 
de cette derni\`ere, la relation:
$$
\Phi_{i}(\mathcal{B}'_{i}) \equiv \mathcal{B}'_{i} B_{i} \pmod{E_{i-1}}.
$$
On en tire enfin la forme triangulaire sup\'erieure par blocs de la matrice
de $\Phi$ dans la base $\mathcal{B}$:
$$
\Phi(\mathcal{B}) = \mathcal{B} 
\begin{pmatrix}
B_{1} & \star  & \star \\ 
0     & \ddots & \star \\
0     & 0      & B_{k} 
\end{pmatrix}.
$$
Conform\'ement aux conventions du paragraphe \ref{subsection:descmat1},
nous identifierons $P_{i}$ ($1 \leq i \leq k$) \`a $(K^{r_{i}},\Phi_{A_{i}})$, 
o\`u $A_{i} := B_{i}^{-1} \in GL_{r_{i}}(K)$. De m\^eme, $P$ s'identifie \`a 
$(K^{n},\Phi_{A_{0}})$ et $M$ \`a $(K^{n},\Phi_{A})$, o\`u 
$n := r_{1} + \cdots + r_{k}$ et:
$$
A_{0} = \begin{pmatrix} 
A_{1} & 0      & 0 \\
0     & \ddots & 0 \\ 
0     & 0      & A_{k} 
\end{pmatrix} \text{~~et~~}
A = \begin{pmatrix} 
A_{1} & \star  & \star \\
0     & \ddots & \star \\ 
0     & 0      & A_{k} 
\end{pmatrix}.
$$
Il faut noter qu'une telle \'ecriture sous-entend implicitement la donn\'ee 
d'une filtration sur $M$ et d'un isomorphisme de $\gr M$ sur $P$. Si de plus 
$M' = (K^{n},\Phi_{A'})$, o\`u $A'$ est de la m\^eme forme que $A$ (\ie\ $M'$
est filtr\'e et muni d'un isomorphisme de $\gr M'$ sur $P$), alors un morphisme 
de $M$ dans $M'$ qui respecte les filtrations (\ie\ qui envoie chaque $M_{i}$ 
dans $M'_{i}$) est d\'ecrit par une matrice de la forme $F$ suivante, et 
l'endomorphisme induit de $P \simeq \gr M \simeq \gr M'$ est d\'ecrit par 
la matrice $F_{0}$ correspondante:
$$
F = \begin{pmatrix} 
F_{1} & \star  & \star \\
0     & \ddots & \star \\ 
0     & 0      & F_{k} 
\end{pmatrix} \text{~~et~~}
F_{0} = \begin{pmatrix} 
F_{1} & 0      & 0     \\
0     & \ddots & 0     \\ 
0     & 0      & F_{k} 
\end{pmatrix}.
$$
En particulier, un morphisme qui induit l'identit\'e sur $P$ (et donc assure
que les modules filtr\'es $M$ et $M'$ appartiennent \`a la m\^eme classe de
$\F(P_{1},\ldots,P_{k})$) est repr\'esent\'e par une matrice de $\G(K)$, o\`u
$\G$ d\'esigne le sous-groupe alg\'ebrique de $GL_{n}(K)$ d\'efini par le
format suivant:
$$
\begin{pmatrix} I_{r_1} & \star  & \star \\
0     & \ddots & \star \\ 0     & 0      & I_{r_k} 
\end{pmatrix}.
$$
Notons $A_{U}$ la matrice triangulaire par blocs dont la composante diagonale
par blocs est $A_{0}$ et donc les blocs sup\'erieurs rectangulaires correspondants
sont les $U_{i,j} \in M_{r_{i},r_{j}}(K)$ ($1 \leq i < j \leq k$); ici, $U$ est
une notation abr\'eg\'ee pour la famille des $U_{i,j}$. Pour tout $F \in \G(K)$,
la matrice $F[A_{U}]$ est de la forme $A_{V}$, pour une famille de
$V_{i,j} \in M_{r_{i},r_{j}}(K)$. On obtient ainsi une op\'eration \`a gauche du
groupe $\G(K)$ sur l'ensemble $\prod\limits_{1 \leq i < j \leq k} M_{r_{i},r_{j}}(K)$.
Les arguments pr\'ec\'edents se r\'esument ainsi:

\begin{prop}
L'application qui, \`a $U$ associe la classe de $(K^{n},\Phi_{A_{U}})$ induit une
bijection entre, d'une part, le quotient de l'ensemble
$\prod\limits_{1 \leq i < j \leq k} M_{r_{i},r_{j}}(K)$ par l'action du groupe
$\G(K)$ et, d'autre part, l'ensemble $\F(P_{1},\ldots,P_{k})$.
\end{prop}

% 2
 
\section{Extensions de modules aux diff\'erences}

Soit $0 \rightarrow M' \rightarrow M \rightarrow M'' \rightarrow 0$
une suite exacte dans $\D_{K,\sigma}-Mod$. Si $M'$ et $M''$ sont des
modules aux diff\'erences, il en est de m\^eme de $M$ (qui est 
$K$-projectif de rang fini parce que la suite est $K$-scind\'ee). 
Le calcul des extensions dans $\DM$ est donc le m\^eme que dans 
$\D_{K,\sigma}-Mod$. On notera simplement $\Ext(M'',M')$ le groupe 
$\Ext_{\D_{K,\sigma}}(M'',M')$ des classes d'extension de $M''$ par $M'$.
D'apr\`es \cite{BAH}, \S 7, ce groupe est muni d'une structure de $C$-module
qui est bien d\'ecrite dans \emph{loc. cit.} 
\footnote{Et, \`a ma connaissance, nulle part ailleurs.}. 
On va expliciter cette structure dans le cas o\`u $M'$ et $M''$ sont 
des modules aux diff\'erences. \\

Soient donc $M = (E,\Phi)$ et $N = (F,\Psi)$. Toute extension
$0 \rightarrow N \rightarrow R \rightarrow M \rightarrow 0$ de $M$
par $N$ donne lieu, par oubli, \`a une suite exacte de $K$-modules
$0 \rightarrow F \rightarrow G \rightarrow E \rightarrow 0$ telle
que, si $R = (G,\Gamma)$, le diagramme suivant commute:
$$
\begin{CD}
0 @>>> F   @>{i}>> G     @>{j}>>      E  @>>> 0     \\
&  &   @VV{\Psi}V      @VV{\Gamma}V      @VV{\Phi}V                 \\
0 @>>> F   @>{i}>> G     @>{j}>>      E @>>> 0
\end{CD}
$$
(On a encore not\'e $i$ et $j$ les applications $K$-lin\'eaires
sous-jacentes aux morphismes \'eponymes.)
Puisque $E$ est projectif, la suite est scind\'ee et l'on peut
d'embl\'ee identifier $G$ au $K$-module $F \times E$, ce que
nous ferons; on note alors $i(y) = (y,0)$ et $j(y,x) = x$.
Les conditions de compatibilit\'e $\Gamma \circ i = i \circ \Psi$
et $\Phi \circ j = j \circ \Gamma$ entrainent alors:
$$
\Gamma(y,x) = \Gamma_{u}(y,x) := \bigl(\Psi(y) + u(x), \Phi(x)\bigr),
\text{~avec~} u \in \Lin_{\sigma}(E,F),
$$
o\`u nous notons $\Lin_{\sigma}(E,F)$ l'ensemble des applications
$\sigma$-lin\'eaires de $E$ dans $F$. (Autrement dit, $u$ est un
morphisme de groupes tel que $u(\lambda x) = \sigma(\lambda) u(x)$.)
Notant de plus $R_{u} := (F \times E,\Gamma_{u})$, qui est un module 
aux diff\'erences naturellement muni d'une structure d'extension 
de $M$ par $N$, on voit que l'on a ainsi d\'efini une application
surjective:
\begin{align*}
\Lin_{\sigma}(E,F) & \rightarrow \Ext(M,N), \\
u & \mapsto \theta_{u} := \text{~classe de~} R_{u}.
\end{align*}
On peut pr\'eciser \`a quelle condition $u,v \in \Lin_{\sigma}(E,F)$
ont m\^eme image $\theta_{u} = \theta_{v}$, \ie\ \`a quelle condition
les extensions $R_{u}$ et $R_{v}$ sont \'equivalentes. Cela se produit
s'il existe un morphisme $\phi: R_{u} \rightarrow R_{v}$ qui induit
l'identit\'e sur $M$ et $N$, autrement dit, une application lin\'eaire
$\phi: F \times E \rightarrow F \times E$ telle que 
$\Gamma_{v} \circ \phi = \phi \circ \Gamma_{u}$ (puisque c'est un 
morphisme de modules aux diff\'erences) et de la forme 
$(x,y) \mapsto \bigl(y + f(x),x\bigr)$ (car elle induit les identit\'es
de $E$ et de $F$). Sous cette derni\`ere forme, la premi\`ere condition
s'\'ecrit:
$$
\forall (y,x) \in F \times E \;,\;
\bigl(\Psi(y + f(x)) + v(x),\Phi(x)\bigr) =
\bigl(\Psi(y) + u(x) + f(\Phi(x)),\Phi(x)\bigr),
$$
c'est-\`a-dire:
$$
u - v = \Psi \circ f - f \circ \Phi.
$$
Remarquons d'ailleurs que, pour tout $f \in \Lin_{K}(E,F)$, l'application
$t_{\Phi,\Psi}(f) := \Psi \circ f - f \circ \Phi$ est $\sigma$-lin\'eaire
de $E$ dans $F$.

\begin{thm}
\label{thm:descriptionextensions}
L'application $u \mapsto \theta_{u}$ de $\Lin_{\sigma}(E,F)$ sur
$\Ext(M,N)$ est fonctorielle en $M$ et en $N$, $C$-lin\'eaire et 
son noyau est l'image de l'application $C$-lin\'eaire:
\begin{align*}
t_{\Phi,\Psi}: \Lin_{K}(E,F) & \rightarrow \Lin_{\sigma}(E,F), \\
f & \mapsto \Psi \circ f - f \circ \Phi. 
\end{align*}
\end{thm}
\Pr

\noindent \underline{Fonctorialit\'e.} \\
Nous ne la prouverons (et ne l'utiliserons) que du c\^ot\'e
covariant, \ie\ en $N$. Nous invoquons \cite{BAH}, \S 7.1 p. 114 
exemple 3 et \S 7.4, p. 119, prop. 4. Soient $\theta$ la classe 
dans $\Ext(M,N)$ de l'extension 
$0 \overset{i}{\rightarrow} N \rightarrow R 
\overset{j}{\rightarrow} M \rightarrow 0$
et $g: N \rightarrow N'$ un morphisme dans $\DM$. Soit
$$
\begin{CD}
0 @>>> N   @>{i}>>  R     @>{j}>>       M  @>>> 0     \\
&  &   @VV{g}V      @VV{h}V      @VV{\Id_{M}}V                 \\
0 @>>> N'  @>{i'}>> R'    @>{j'}>>      M  @>>> 0
\end{CD}
$$
un diagramme commutatif de suites exactes.
Si $\theta'$ est la classe dans $\Ext(M,N)$ de l'extension 
$0 \overset{i'}{\rightarrow} N' \rightarrow R'
\overset{j'}{\rightarrow} M \rightarrow 0$,
alors: 
$$
\Ext(\Id_{M},g)(\theta) = g \circ \theta = 
\theta' \circ \Id_{M} = \theta'.
$$
On peut prendre par exemple:
$$
R' := R \oplus_{N} N' = 
\dfrac{R \times N'}{\{\bigl(i(n),-g(n)\bigr) \tq n \in N\}},
$$
avec pour $i',j'$ les fl\`eches \'evidentes. Prenant pour $R$
l'extension $R_{u}$ et reprenant les notations ant\'erieures
pour $N,M,R$, et notant de plus $N' = (F',\Psi')$, avec la
condition de compatibilit\'e $\Psi' \circ g = g \circ \Psi$,
on voit que le $K$-module sous-jacent \`a $R \oplus_{N} N'$
est:
$$
G' := \dfrac{F \times E \times F'}{\{\bigl(y,0,-g(y)\bigr) \tq y \in F\}},
$$
muni de l'automorphisme semi-lin\'eaire induit par l'application
$\Gamma_{u} \times \Psi'$ de $F \times E \times F'$ dans lui-m\^eme
(celle-ci laisse bien stable le d\'enominateur). \\
L'application $(y,x,y') \mapsto \bigl(y' + g(y),x\bigr)$ de 
$F \times E \times F'$ dans $F' \times E$ induit un isomorphisme
de $G'$ sur $F' \times E$ et l'automorphisme semi-lin\'eaire
induit sur $G'$ est 
$(y',x) \mapsto \bigl(\Psi'(y') + g\bigl(u(x)\bigr), \Phi(x)\bigr)$,
c'est-\`a-dire $\Gamma_{gu}$, d'o\`u l'on tire que $R' = R_{gu}$.
Les fl\`eches $i'$ et $j'$ sont d\'etermin\'ees comme suit:
$i'(y')$ est la classe de $(0,y')$ dans $G'$, soit, avec 
l'identification ci-dessus, $i'(y') = (y',0)$; et $j'(y',x)$ est
l'image d'un ant\'ec\'edent quelconque, par exemple la classe
de $(0,x,y')$, image qui vaut $j(0,x) = x$. On a donc \'etabli
que la classe de l'extension $R_{u}$ par $\Ext(\Id_{M},g)$ est
$R_{gu}$, ce qui est la fonctorialit\'e annonc\'ee. Elle se
traduit par la commutativit\'e du diagramme suivant:
$$
\begin{CD}
\Lin_{\sigma}(E,F)  @>>> \Ext(M,N) \\
@VV{\Lin_{\sigma}(\Id_{M},g)}V   @VV{\Ext(\Id_{M},g)}V  \\
\Lin_{\sigma}(E,F') @>>> \Ext(M,N')
\end{CD}
$$

\noindent \underline{Lin\'earit\'e.} \\
L'application $t_{\Phi,\Psi}$ va bien de $\Lin_{K}(E,F)$ dans
$\Lin_{\sigma}(E,F)$ d'apr\`es la remarque qui pr\'ec\`ede
imm\'ediatement du th\'eor\`eme. \\
\emph{Addition.}
La r\'ef\'erence est ici \cite{BAH},\S 7.6, rem. 2 p. 124.
\`A partir des extensions
$0 \rightarrow N \overset{i}{\rightarrow} R 
\overset{p}{\rightarrow} M \rightarrow 0$ and
$0 \rightarrow N \overset{i'}{\rightarrow} R'
\overset{p'}{\rightarrow} M \rightarrow 0$ de classes 
$\theta,\theta' \in \Ext^{1}(M,N)$, on calcule $\theta + \theta'$
comme classe de l'extension
$0 \rightarrow N \overset{i''}{\rightarrow} R''
\overset{p''}{\rightarrow} M \rightarrow 0$, o\`u:
$$
R'' := \dfrac{\{(z,z') \in R \times R' \tq p(z) = p'(z')\}}
{\{(-i(y),i'(y)) \tq y \in N\}},
$$
et $i''(y)$ est la classe de $(0,i'(y))$, \ie\ la m\^eme que la classe de
$(i(y),0)$; et $p''$ envoie la classe de $(z,z')$ sur $p(z) = p'(z')$.
Si l'on prend $R = R_{u}$ et $R' = R_{u'}$, le num\'erateur de $R''$ 
s'identifie \`a $F \times F \times E$ muni de l'automorphisme 
semi-lin\'eaire $(y,y',x) \mapsto (\Psi(y) + u(x),\Psi(y') + u'(x),\Phi(x))$. 
Le d\'enominateur s'identifie au sous-espace $\{(-y,y,0) \tq y \in F\}$
muni de l'application induite. Le quotient s'identifie \`a
$0 \times F \times E$, via l'application $(y,y',x) \mapsto (0,y'',x)$,
o\`u $y'' := y' + y$ muni de l'automorphisme semi-lin\'eaire $\Phi''$ qui
envoie $(0,y'',x)$ sur
$$
(0,\Psi(y') + u'(x) + \Psi(y) + u(x),\Phi(x)) =
(0,\Psi(y'') + (u + u')(x),\Phi(x)).
$$
On reconnait bien $R_{u+u'}$. \\
\emph{Multiplication externe.}
La r\'ef\'erence est ici \cite{BAH},\S 7.6, prop. 4 p. 119. 
Soit $\lambda \in C$. On applique la proposition indiqu\'ee
au diagramme commutatif de suites exactes:
$$
\begin{CD}
0 @>>> N   @>>>  R_{u}          @>>>       M  @>>> 0     \\
&  &  @VV{\times \lambda}V  @VV{(\times \lambda,\Id_{M})}V @VV{\Id_{M}}V \\
0 @>>> N  @>>> R_{\lambda u}    @>>>       M  @>>> 0
\end{CD}
$$
Si $\theta$, $\theta'$ sont les classe dans $\Ext(M,N)$ des deux
extensions, on d\'eduit de \emph{loc. cit.} que:
$$
\theta' \circ \Id_{M} = (\times \lambda) \circ \theta 
\Longrightarrow \theta' = \lambda \theta.
$$
La classe de l'extension $R_{\lambda u}$ est donc bien \'egale
au produit par $\lambda$ de la classe de l'extension $R_{u}$. \\

\noindent \underline{Exactitude.} \\
Elle d\'ecoule imm\'ediatement du calcul qui pr\'ec\`ede l'\'enonc\'e
du th\'eor\`eme.
\Endpr

% 2.1

\subsection{Le complexe des solutions}
\label{subsection:complexedessolutions}

\begin{defn}
On appelle \emph{complexe des solutions de $M$ dans $N$}
le complexe de $C$-modules:
\begin{align*}
t_{\Phi,\Psi}: \Lin_{K}(E,F) & \rightarrow \Lin_{\sigma}(E,F), \\
f & \mapsto \Psi \circ f - f \circ \Phi. 
\end{align*}
concentr\'e aux degr\'es $0$ et $1$.
\end{defn}

Il est bien clair que ce sont bien des $C$-modules et que 
l'application est bien d\'efinie (elle va dans son but !)
et $C$-lin\'eaire. 

\begin{cor}
Le complexe des solutions a pour homologie $H^{0} = \Hom(M,N)$ 
et $H^{1} = \Ext(M,N)$ et ceci, fonctoriellement.
\end{cor}
\Pr 
Pour $H^{1}$, c'est le th\'eor\`eme. Pour $H^{0}$, on voit que
le noyau de $t_{\Phi,\Psi}$ est le $C$-module 
$\{f \in \Lin_{K}(E,F) \tq \Psi \circ f = f \circ \Phi\}$,
c'est-\`a-dire $\Hom(M,N)$; et la fonctorialit\'e est ici \'evidente.
\Endpr

\begin{cor}
\label{cor:longuesuite}
De la suite exacte 
$0 \rightarrow N' \rightarrow N \rightarrow N'' \rightarrow 0$,
on d\'eduit la ``longue suite exacte de cohomologie'':
$$
0 \rightarrow \Hom(M,N') \rightarrow \Hom(M,N) \rightarrow \Hom(M,N'') 
\rightarrow \Ext(M,N') \rightarrow \Ext(M,N) \rightarrow \Ext(M,N'')
\rightarrow 0.
$$
\end{cor}
\Pr
On conserve les notations pr\'ec\'edentes (que l'on adapte de plus
\`a $N''$). La suite exacte de $K$-modules (projectifs)
$0 \rightarrow F' \rightarrow F \rightarrow F'' \rightarrow 0$
\'etant scind\'ee, les deux lignes du diagramme commutatif:
$$
\begin{CD}
0 @>>> \Lin_{K}(E,F') @>>> \Lin_{K}(E,F)  @>>> \Lin_{K}(E,F'')  @>>> 0     \\
&  &  @VV{t_{\Phi,\Psi'}}V @VV{t_{\Phi,\Psi}}V @VV{t_{\Phi,\Psi''}}V  \\
0 @>>> \Lin_{\sigma}(E,F') @>>> \Lin_{\sigma}(E,F)  @>>> 
\Lin_{\sigma}(E,F'')  @>>> 0 
\end{CD}
$$
sont exactes, et il suffit alors d'invoquer le lemme du serpent.
\Endpr

% 2.2

\subsection{Description matricielle des extensions de modules 
aux diff\'erences}
\label{subsection:descmat3}

On suppose ici $E$ et $F$ libres sur $K$, et l'on op\`ere les 
identifications correspondantes $M = (K^{m},\Phi_{A})$, $A \in GL_{m}(K)$
et $N = (K^{n},\Phi_{B})$, $B \in GL_{n}(K)$. Une extension de $N$ par $M$ 
est alors de la forme $R = (K^{m+n},\Phi_{C})$, o\`u
$C = \begin{pmatrix} A & U \\ 0_{n,m} & B \end{pmatrix}$ pour une certaine
matrice rectangulaire $U \in M_{m,n}(K)$; nous noterons $C = C_{U}$.
L'injection $M \rightarrow R$ et la projection $R \rightarrow N$ ont
respectivement pour matrice $\begin{pmatrix} I_{m} \\ 0_{n,m} \end{pmatrix}$
et $\begin{pmatrix} 0_{n,m} & I_{n}\end{pmatrix}$. L'extension ainsi 
d\'efinie sera not\'ee $R_{U}$. \\

Un morphisme d'extensions $R_{U} \rightarrow R_{V}$ est une matrice de la
forme $F = \begin{pmatrix} I_{m} & X \\ 0_{n,m} & I_{n} \end{pmatrix}$ pour
une certaine matrice rectangulaire $X \in M_{m,n}(K)$. La condition de
compatibilit\'e avec les automorphismes semi-lin\'eaires s'\'ecrit:
$$
(\sigma F) C_{U} = C_{V} F \Longleftrightarrow U + (\sigma X) B = A X + V
\Longleftrightarrow V - U = (\sigma X) B - A X.
$$

\begin{cor}
\label{cor:descmatrext}
Le $C$-module $\Ext^{1}(N,M)$ s'identifie ainsi avec le conoyau de 
l'endomorphisme $X \mapsto (\sigma X) B - A X$ de $M_{m,n}(K)$.
\end{cor}
\Pr
La construction ci-dessus fournit une bijection, mais il r\'esulte
du th\'eor\`eme \ref {thm:descriptionextensions} qu'il s'agit bien 
d'un isomorphisme.
\Endpr

% 3

\section{Extension des scalaires}

On veut consid\'erer $\F(P_{1},\ldots,P_{k})$ comme un sch\'ema
sur $C$, donc comme un foncteur (repr\'esentable)
$C' \leadsto \F(C' \otimes_{C} P_{1},\ldots,C' \otimes_{C} P_{k})$
des $C$-alg\`ebres commutatives vers les ensembles. Pour cela,
on va \'etendre ce qui pr\'ec\`ede en une situation ``relative''. \\

Soit $C'$ une $C$-alg\`ebre commutative. On note: 
$$
K' := C' \otimes_{C} K \text{~et~} \sigma' := 1 \otimes_{C} \sigma.
$$
Donc $K'$ est une $C'$-alg\`ebre commutative et $\sigma'$ est
un automorphisme de cette $C'$-alg\`ebre. On a alors:
$$
\D_{K',\sigma'} := K'<\sigma',{\sigma'}^{-1}> = 
K' \otimes_{K} \D_{K,\sigma} = C' \otimes_{C} \D_{K,\sigma}.
$$
Ces \'egalit\'es signifient: isomorphismes naturels (fonctoriels). \\

Pour tout module aux $q$-diff\'erences $M = (E,\Phi)$ sur $(K,\sigma)$, 
on obtient un module aux $q$-diff\'erences $M' = (E',\Phi')$ sur 
$(K',\sigma')$ en posant:
$$
E' = K' \otimes_{K} E = C' \otimes_{C} E \text{~et~}
\Phi' = \sigma' \otimes_{K} \Phi = 1 \otimes_{C} \Phi.
$$
(C'est bien un $\D_{K',\sigma'}$-module \`a gauche qui est projectif
de rang fini sur $K'$.) Nous le noterons $M' = C' \otimes_{k} M$
pour faire ressortir sa d\'ependance en $C'$. La proposition suivante
est l'outil pour traiter le cas o\`u $k = 2$.

\begin{prop}
\label{prop:supermajik}
Soient $M,N$ deux modules aux $q$-diff\'erences sur $(K,\sigma)$.
On a un isomorphisme fonctoriel de $C'$-modules:
$$
\Ext_{\D_{K',\sigma'}}(C' \otimes_{C} M,C' \otimes_{C} N) 
\simeq C' \otimes_{C} \Ext_{\D_{K,\sigma}}(M,N),
$$
et un \'epimorphisme fonctoriel de $C'$-modules:
$$
C' \otimes_{C} \Hom_{\D_{K,\sigma}}(M,N) \rightarrow
\Hom_{\D_{K',\sigma'}}(C' \otimes_{C} M,C' \otimes_{C} N).
$$
\end{prop}
\Pr
On notera $M' = C' \otimes_{C} M$, $E' = K' \otimes_{K} E$ etc.
Les $K$-modules $E,F$ \'etant projectifs de rang fini, on a des
isomorphismes naturels:
$$
C' \otimes_{C} \Lin_{K}(E,F) = \Lin_{K'}(E',F') \text{~et~}
C' \otimes_{C} \Lin_{\sigma}(E,F) = \Lin_{\sigma'}(E',F').
$$
(C'est imm\'ediat si $E$ et $F$ sont libres et le cas g\'en\'eral 
s'en d\'eduit.) En tensorisant la suite exacte (fonctorielle):
$$
0 \rightarrow \Hom_{\D_{K,\sigma}}(M,N) \rightarrow \Lin_{K}(E,F) 
\rightarrow \Lin_{\sigma}(E,F) \rightarrow \Ext_{\D_{K,\sigma}}(M,N)
\rightarrow 0,
$$
on obtient la suite exacte:
$$
C' \otimes_{C} \Hom_{\D_{K,\sigma}}(M,N) \rightarrow 
C' \otimes_{C} \Lin_{K}(E,F) \rightarrow C' \otimes_{C} \Lin_{\sigma}(E,F) 
\rightarrow C' \otimes_{C} \Ext_{\D_{K,\sigma}}(M,N) \rightarrow 0.
$$
Les deux conclusions viennent alors par comparaison avec la suite exacte:
$$
0 \rightarrow \Hom_{\D_{K',\sigma'}}(M',N') \rightarrow \Lin_{K'}(E',F') 
\rightarrow \Lin_{\sigma'}(E',F') \rightarrow \Ext_{\D_{K',\sigma'}}(M',N')
\rightarrow 0.
$$
\Endpr

\begin{prop}
Soit $0 = M_{0} \subset M_{1} \subset \cdots \subset M_{k} = M$
une $k$-filtration de gradu\'e $P_{1} \oplus \cdots \oplus P_{k}$.
Alors, notant $M'_{i} := C' \otimes_{C} M_{i}$ et 
$P'_{i} := C' \otimes_{C} P_{i}$, on obtient une $k$-filtration 
$0 = M'_{0} \subset M'_{1} \subset \cdots \subset M'_{k} = M'$
de gradu\'e $P'_{1} \oplus \cdots \oplus P'_{k}$.
\end{prop}
\Pr
Les $P_{i}$ \'etant projectifs en tant que $K$-modules, les suites
exactes 
$0 \rightarrow M_{i-1} \rightarrow M_{i} \rightarrow P_{i} \rightarrow 0$
sont $K$-scind\'ees, donc donnent lieu par changement de base
$K' \otimes_{K}$ aux suites exactes
$0 \rightarrow M'_{i-1} \rightarrow M'_{i} \rightarrow P'_{i} \rightarrow 0$.
\Endpr

Si $(\underline{M},\underline{u})$ d\'esigne le couple form\'e
de l'objet $k$-filtr\'e ci-dessus et d'un isomorphisme pr\'ecis\'e
de $\gr M$ sur $P_{1} \oplus \cdots \oplus P_{k}$, on notera 
$(C' \otimes_{C} \underline{M},1 \otimes_{C} \underline{u})$ 
le couple analogue d\'eduit de la proposition.

\begin{defn}
On d\'efinit comme suit un foncteur $F$ de la cat\'egorie 
des $C$-alg\`ebres commutatives dans la cat\'egorie des ensembles.
Pour toute $C$-alg\`ebre commutative $C'$, on pose:
$$
F(C') := \F(C' \otimes_{C} P_{1},\ldots,C' \otimes_{C} P_{k}).
$$
Pour tout morphisme $C' \rightarrow C''$ de $C$-alg\`ebres commutatives,
l'application $F(C') \rightarrow F(C'')$ est donn\'ee par:
$$
\text{classe de~} (\underline{M'},\underline{u'}) \mapsto
\text{classe de~} 
(C'' \otimes_{k'} \underline{M'},1 \otimes_{C} \underline{u'}).
$$
\end{defn}

L'ensemble $F(C')$ est bien d\'efini d'apr\`es les constructions 
du d\'ebut de cette section. L'application $F(C') \rightarrow F(C'')$
est bien d\'efinie au niveau des couples en vertu de la proposition
et elle passe au quotient (v\'erification laiss\'ee au lecteur).
Enfin, on a bien un foncteur (l'image d'un morphisme compos\'e est
le compos\'e des images) en vertu de la r\`egle de contraction des
produits tensoriels, qui donne ici:
$$
C''' \otimes_{C''} (C'' \otimes_{C'} \underline{M'}) =
C''' \otimes_{C'} \underline{M'}.
$$

% 4

\section{Notre espace de modules}

Pour simplifier (et parce que c'est correct !), dans ce qui suit,
au lieu de dire ``le foncteur $F$ est est repr\'esentable par un 
espace affine sur $C$ (de dimension $d$)'', on dira ``le foncteur
$F$ est un espace affine sur $C$ (de dimension $d$)''

\begin{thm}
On suppose que, pour $1 \leq i < j \leq k$, on a $\Hom(P_{j},P_{i}) = 0$
et que le $C$-module $\Ext(P_{j},P_{i})$ est libre de rang fini
$\delta_{i,j}$. Alors le foncteur 
$C' \leadsto F(C') := \F(C' \otimes_{C} P_{1},\ldots,C' \otimes_{C} P_{k})$
est un espace affine sur $C$ de dimension
$\sum\limits_{1 \leq i < j \leq k} \delta_{i,j}$.
\end{thm}
\Pr
Lorsque $k = 1$, c'est trivial. Lorsque $k = 2$, notant $V$ le 
$C$-module libre de rang fini $\Ext(P_{2},P_{1})$, il s'agit
(en vertu de la proposition \ref{prop:supermajik}) du foncteur 
$C' \leadsto C' \otimes_{C} V$ qui est repr\'esent\'e par l'alg\`ebre 
sym\'etrique du dual de $V$, laquelle est une alg\`ebre de polyn\^omes 
sur $C$. Au del\`a, on va raisonner par r\'ecurrence sur $k$ en 
invoquant le lemme 2.5.3, p. 139 de \cite{BV}:
\begin{lem}
Soit $u: F \rightarrow G$ une transformation naturelle entre deux
foncteurs des $C$-alg\`ebres commutatives vers les ensembles.
On suppose que $G$ est un espace affine sur $C$ et que, pour 
toute $C$-alg\`ebre commutative $C'$ et tout $b \in G(C')$, 
la ``fibre de $u$ en $b$'', qui est le foncteur des $C'$-alg\`ebres 
commutatives vers les ensembles
$$
C'' \leadsto u_{C''}^{-1}\bigl(G(C' \rightarrow C'')(b)\bigr)
$$
est un espace affine sur $C'$. Alors $F$ est un espace affine sur $C$.
\end{lem}
Dans \emph{loc. cit.}, ce th\'eor\`eme est prouv\'e pour $C = \C$,
mais l'argument est manifestement valable pour tout anneau commutatif.
En voici le squelette. On choisit $B = C[T_{1},\ldots,T_{d}]$ qui
repr\'esente (ou dont le spectre repr\'esente) $G$. On prend pour $b$
l'identit\'e de $G(B) = \Hom(B,B)$ (``point g\'en\'eral''); la fibre
est repr\'esent\'ee par $B[S_{1},\ldots,S_{e}]$. On prouve alors
que  $C[T_{1},\ldots,T_{d},S_{1},\ldots,S_{e}]$ repr\'esente $F$.
Cela donne en passant un calcul de la dimension de l'espace affine
$F$ comme somme des dimensions de $G$ et de la ``fibre g\'en\'erale''. 
Dans notre application, toutes les fibres auront m\^eme dimension. 
Avant de poursuivre la preuve du th\'eor\`eme, d\'emontrons une 
proposition auxiliaire.
\begin{prop}
Soit $C'$ une $C$-alg\`ebre commutative et soit $M'$ un module 
aux diff\'erences sur $K' := C' \otimes_{C} K$ admettant une 
$(k-1)$-filtration:
$0 = M'_{0} \subset M'_{1} \subset \cdots \subset M'_{k-1} = M'$
telle que $\gr M' \simeq P'_{1} \oplus \cdots \oplus P'_{k-1}$
(comme d'habitude, $P'_{i} := C' \otimes_{C} P_{i}$). Alors le
foncteur en $C'$-alg\`ebres commutatives
$C'' \leadsto \Ext(C'' \otimes_{C} P_{k}, C'' \otimes_{C'} M')$
est un espace affine sur $C'$ de dimension
$\sum\limits_{1 \leq i \leq k} \delta_{i,k}$.
\end{prop}
\Pr
D'apr\`es la proposition \ref{prop:supermajik}, il s'agit du foncteur
$C'' \leadsto C'' \otimes_{C'} \Ext(P'_{k},M')$. De chaque suite exacte
$0 \rightarrow M'_{i-1} \rightarrow M'_{i} \rightarrow P'_{i} \rightarrow 0$
on d\'eduit la longue suite exacte de cohomologie d\'ecrite au corollaire
\ref{cor:longuesuite}; mais, de la proposition \ref{prop:supermajik}, on 
d\'eduit que, pour toute $C$-alg\`ebre commutative $C'$, on a 
$\Hom(C' \otimes_{C} P_{j},C' \otimes_{C} P_{i}) = 0$ et que 
le $C'$-module $\Ext(C' \otimes_{C} P_{j},C' \otimes_{C} P_{i})$ 
est libre de rang fini $\delta_{i,j}$. Vues les \'egalit\'es 
$\Hom(P'_{j},P'_{i}) = 0$, la longue suite exacte prend ici 
la forme raccourcie:
$$
0 \rightarrow \Ext(P'_{k},M'_{i-1}) \rightarrow \Ext(P'_{k},M'_{i}) 
\rightarrow \Ext(P'_{k}, P'_{i}) \rightarrow 0,
$$
et, pour $i = 1,\ldots,k-1$, ces suites sont scind\'ees, le membre
droit \'etant libre. On a donc finalement:
$$
\Ext(P'_{k},M') \simeq \bigoplus_{1 \leq i \leq k} \Ext(P'_{k}, P'_{i}),
$$
qui est libre de rang $\sum\limits_{1 \leq i \leq k} \delta_{i,k}$.
Comme dans le cas $k = 2$ (qui est d'ailleurs un cas particulier de
cette proposition), le foncteur indiqu\'e est repr\'esentable
par l'alg\`ebre sym\'etrique du dual de ce module.
\Endpr
Terminons la preuve du th\'eor\`eme. Outre le foncteur $F(C')$,
on consid\`ere le foncteur
$C' \leadsto G(C') := \F(C' \otimes_{C} P_{1},\ldots,C' \otimes_{C} P_{k-1})$,
dont on peut supposer que c'est un espace affine sur $C$ de dimension 
$\sum\limits_{1 \leq i < j \leq k-1} \delta_{i,j}$ (hypoth\`ese
de r\'ecurrence). La transformation naturelle de $F$ dans $G$ est celle
que l'on a d\'ecrite \`a la fin de l'introduction, page \pageref{methode}.
Un \'el\'ement $b \in G(C')$ est la classe d'un couple
$(\underline{M'},\underline{u'})$, objet $(k-1)$-filtr\'e sur $C'$,
et la fibre correspondante est le foncteur \'etudi\'e dans la proposition
auxiliaire ci-dessus. Le lemme de Babbitt et Varadarajan permet
alors de conclure.
\Endpr

\Rem
L'isomorphisme 
$\Ext(P'_{k},M') \simeq 
\bigoplus\limits_{1 \leq i \leq k} \Ext(P'_{k}, P'_{i})$
n'est pas fonctoriel. Cependant, on doit pouvoir tirer un peu plus
de la preuve ci-dessus: 
\begin{itemize}
\item{La transformation naturelle de $F$ dans $G$ identifie l'espace
affine $F$ au produit de l'espace affine $G$ par l'espace affine
$C' \leadsto \bigoplus\limits_{1 \leq i \leq k} \Ext(P'_{k}, P'_{i})$.}
\item{Si l'on a des coordonn\'ees sur chaque $\Ext(P_{j}, P_{i})$,
on a des coordonn\'ees sur $F$.}
\end{itemize}
Ces am\'eliorations sont exploit\'ees dans \cite{RSZ}.

%%%%%%%%%%%%%%%%%%%%%%%%%%%%%%%%%%%%%%%%%%%%%%%%%%%%%%%%%%%%%%%%%%%%%%%%%%%%%

% Biblio

\end{document}